\documentclass{article}
\usepackage{amsfonts}

\usepackage{amsmath}

\newtheorem{theorem}{Theorem}

\newtheorem{corollary}[theorem]{Corollary}

\newtheorem{lemma}[theorem]{Lemma}

\newtheorem{proposition}[theorem]{Proposition}
\newtheorem{remark}[theorem]{Remark}

\newenvironment{proof}[1][Proof]{\noindent\textbf{#1.} }{\ \rule{0.5em}{0.5em}}

\input{tcilatex}

\begin{document}

\title{Linearizability of $d$-webs, $d\geq 4,$ on two-dimensional manifolds}
\author{Maks A. Akivis, Vladislav V. Goldberg, Valentin V. Lychagin}
\maketitle

\begin{abstract}
We find $d - 2$ relative differential invariants for a $d$-web, $d\geq 4,$
on a two-dimensional manifold and prove that their vanishing is necessary
and sufficient for a $d$-web to be linearizable. If one writes the above
invariants in terms of web functions $f(x,y)$ and $%
g_{4}(x,y),...,g_{d}(x,y), $ then necessary and sufficient
conditions for the linearizabilty of a $d$-web are two PDEs of the
fourth order with respect to $f$ and $g_{4}$, and $d - 4$ PDEs of
the second order with respect to $f$ and $g_{4},...,g_{d}$. For $d
= 4,$ this result confirms Blaschke's conjecture on the nature of
conditions for the linearizabilty of a $4$-web. We also give
Mathematica codes for testing $4$- and $d$-webs ($d
> 4$) for linearizability and examples of their usage.
\end{abstract}

\section{\protect\bigskip Introduction}

Let $W_{d}$ be a $d$-web given by $d$ one-parameter foliations of curves on
a two-dimensional manifold $M^{2}$. The web $W_{d}$ is linearizable
(rectifiable) if it is equivalent to a linear $d$-web, i.e., to a $d$-web
formed by $d$ one-parameter foliations of straight lines on a projective
plane.

The problem of the linearizability of webs was posed by Blaschke in the
1920s (see, for example, his book \cite{B 55}, \S 17 and \S 42) who claimed
that it is hopeless to find such a criterion because of the complexity of
calculations involving high order jets. \ Blaschke in \cite{B 55} (\S\ 42)
formulated the problems of finding conditions for the linearizability of $3$%
-webs (\S\ 17) and $4$-webs (\S\ 42) given on $M^{2}.$ \ Comparing \ the
numbers of\ absolute invariants for a general $3$-web $W_{3}$ (a general $4$%
-web $W_{4}$) and a linear $3$-web (a linear $4$-web), Blaschke made the
conjectures that conditions of linearizability for a $3$-web $W_{3}$ should
consist of four relations for the $9$th order web invariants ($4$ PDE of $9$%
th order) and those for a $4$-web $W_{4}$ should consist of two relations
for the $4$th order web invariants ($2$ PDE of $4$th order) .

A criterion for linearizability is very important in web geometry and in its
applications. It is also important in applications to nomography (see \cite%
{B 55}, \S 17 and \cite{BB 38}, \S 18).

A new approach for finding conditions of linearizability for webs on the
plane has been proposed by Akivis (1973) in his talk at the Seminar on
Classical Differential Geometry in Moscow State University. Goldberg \cite{G
89} implemented this approach for $3$-webs. The goals of the authors of the
paper \cite{GMS 01} were to find linearizability conditions for a $3$-web $%
W_{3}$ and to improve Bol's and Boruvka's result related to the Gronwall
conjecture. For the formulation of\ the Gronwall conjecture, the statement
of the results of Bol and Boruvka and references to their works see \cite{B
55}, \S 17.

In this paper we use Akivis' approach to establish a criterion of
linearizability for $d$-webs, $d\geq 4$. \ The results of the present paper
do not rely on the results or methods of the paper \cite{GMS 01} mentioned
above.\ We prove that the Blaschke conjecture was correct: a $4$-web $W_{4}$
is linearizable if and only its two $4$th order invariants vanish. In terms
of the invariants defining the geometry of a $4$-web $W_{4},$ the vanishing
of these two invariants means that the covariant derivatives $K_{1}$ and $%
K_{2}$\ of the web curvature $K$\ are expressed in terms of the curvature $K$
itself, the basic web invariant $a$ and its covariant derivatives up to the $%
3$rd order. We find explicit expressions for these invariants in terms of
symmetrized covariant derivatives. Note that expressions for these
invariants in terms of web functions contain $262$ terms each. After this
paper was submitted, one of the authors used the conditions of
linearizability described above to check whether numerous known classes of $4
$-webs are linearizable (see \cite{G 04}).

The results obtained in this paper give a complete solution of the
linearizability problem for $d$-webs, $d\geq 4$, and provide tests for
establishing linearizability of such webs. In particular, for $4$-webs $%
W_{4},$ our results provide a complete solution of the longstanding problem
posed by Blaschke (see, for example \cite{B 55}, \S 42).

We also investigate the linearizability of $d$-webs $W_{d}$ for $d\geq 5.$
In this case the linearizability conditions involve $d-2$ differential
invariants. Two of them have order $4$ and the rest are of order $2$.

All computations in this paper were done manually, and the more routine ones
(for example, equations (13), (14), 15) and the formulas for $K_{1}$ and $%
K_{2}$ in Section 2.3.4) were checked by Mathematica package. At the end of
the paper, we provide the Mathematica codes for testing $4$- and $d$-webs, $%
d>4,$ for linearizability and examples of their usage. The material in
Section $4$ (tests and examples) essentially relies on using Mathematica.

Note that a different approach to the linearizability problem for webs $W_{d}
$ for $d\geq 4$ was used by H\`{e}naut in \cite{H 93}. However, H\`{e}naut
did not find conditions in the form suggested by Blaschke. His conditions do
not contain web invariants.

\section{Basics Constructions}

We recall main constructions for $3$-webs on $2$-dimensional manifolds (see,
for example, \cite{BB 38} or \cite{B 55} , or \cite{G 89}) in a form
suitable for us.

Let $M^{2}$ be a $2$-dimensional manifold, and suppose that a $3$-web $W_{3}$
is given by three differential $1$-forms $\omega _{1},\omega _{2},$ and $%
\omega _{3}$ such that any two of them are linearly independent.

\begin{proposition}
The forms $\omega _{1},\omega _{2},$ and $\omega _{3}$ can be normalized in
such a way that the normalization condition
\begin{equation}
\omega _{1}+\omega _{2}+\omega _{3}=0  \label{normalization equation}
\end{equation}%
holds.
\end{proposition}

\begin{proof}
In fact, if we take the forms $\omega _{1}$ and $\omega _{2}$ as cobasis
forms of $M^{2}$, then the form $\omega _{3}$ is a linear combination of the
forms $\omega _{1}$ and $\omega _{2}:$
\begin{equation*}
\omega _{3}=\alpha \omega _{1}+\beta \omega _{2\,},
\end{equation*}
where $\alpha ,\beta \neq 0.$

After the substitution
\begin{equation*}
\omega _{1}\rightarrow \frac{1}{\alpha }\omega _{1},\;\omega _{2}\rightarrow
\frac{1}{\beta }\omega _{2},\;\omega _{3}\rightarrow -\omega _{3}
\end{equation*}
the above equation becomes (\ref{normalization equation}).
\end{proof}

It is easy to see that any two of such normalized triplets $\omega
_{1},\omega _{2},\omega _{3}$ and $\omega _{1}^{s},\omega _{2}^{s},\omega
_{3}^{s}$ determine the same $3$-web $W_{3}$ if and only if
\begin{equation}
\omega _{1}^{s}=s\omega _{1},\ \omega _{2}^{s}=s\omega _{2},\ \omega
_{3}^{s}=s\omega _{3}  \label{perenormirovka}
\end{equation}
for a non-zero smooth function $s.$

We will investigate local properties of $W_{3}.$ Thus we can assume that $%
M^{2}$ is a simply connected domain of $\mathbb{R}^{2}$, and therefore there
exists a smooth function $f$ such that $\omega _{3}$ is proportional to $df,$
that is, $\omega _{3}\wedge df=0.$ The function $f$ is called a \emph{web
function. }Note that this function is defined up to renormalization $%
f\longmapsto F\left( f\right) .$

We choose such a representation of $W$ that
\begin{equation}
\omega _{3}=df.  \label{normalized condition2}
\end{equation}

Similarly we find smooth functions $x$ and $y$ for forms $\omega _{1}$ and $%
\omega _{2}$ such that%
\begin{equation*}
\omega _{1}=adx,\ \omega _{2}=bdy
\end{equation*}%
for some smooth functions $a$ and $b.$

Moreover, functions $x$ and $y$ are independent and therefore can be viewed
as (local) coordinates. In these coordinates the normalization condition
gives
\begin{equation*}
\omega _{1}=-f_{x}dx,\ \omega _{2}=-f_{y}dy,\ \omega _{3}=df.
\end{equation*}%
Let the vector fields $\partial _{1}$ and $\partial _{2}$ form the basis
dual to the cobasis $\omega _{1},\omega _{2},$\ i.e.,$\ \ \ \omega
_{i}\left( \partial _{j}\right) =\delta _{ij}$ for $i,j=1,2.$

Then
\begin{equation*}
\partial _{1}=-\frac{1}{f_{x}}\frac{\partial }{\partial x},\ \ \ \partial
_{2}=-\frac{1}{f_{y}}\frac{\partial }{\partial y}
\end{equation*}
and
\begin{equation}
dv=\partial _{1}(v)~\omega _{1}+\partial _{2}\left( v\right) ~\omega _{2}
\label{differential}
\end{equation}
for any smooth function $v.$

\subsection{Structure Equations}

From now on we shall assume that a $3$-web $W_{3}$ is given by differential $%
1$-forms $\omega _{1},\omega _{2},$ and $\omega _{3}$ normalized by
conditions (\ref{normalization equation}) and (\ref{normalized condition2}).

Since on\ a two-dimensional manifold the exterior differentials $d\omega
_{1} $ and $d\omega _{2}$ as $2$-forms differ from the $2$-form $\omega
_{1}\wedge \omega _{2}$ \ only by factors, we get $d\omega _{1}=h_{1}\
\omega _{1}\wedge \omega _{2}$ and $d\omega _{2}=h_{2}\ \omega _{2}\wedge
\omega _{1}$ for some functions $h_{1}$ and $h_{2}.$

By $d\omega _{3}=0,$ one gets $h_{1}=h_{2}.$ Denote this function by $H.$
Then $d\omega _{1}=H\omega _{1}\wedge \omega _{2}\ \ $and $d\omega
_{2}=H\omega _{2}\wedge \omega _{1}$ or
\begin{equation}
d\omega _{1}=\omega _{1}\wedge \gamma ,\ \ d\omega _{2}=\omega _{2}\wedge
\gamma ,  \label{structure equations}
\end{equation}%
where
\begin{equation}
\gamma =-H\omega _{3}.  \label{Gamma}
\end{equation}%
We call relations (\ref{structure equations}) the \emph{first structure
equations} of the $3$-web $W_{3}.$ In terms of the web function $f$, one has
\begin{equation*}
\gamma =-\frac{f_{xy}}{f_{x}f_{y}}\omega _{3}
\end{equation*}%
and
\begin{equation*}
H=\frac{f_{xy}}{f_{x}f_{y}}.
\end{equation*}%
If we change the representative according to (\ref{perenormirovka}), then
the first structure equations take the form
\begin{equation*}
d\omega _{p}^{s}=\omega _{p}^{s}\wedge \gamma ^{s},\ p=1,2,3,
\end{equation*}%
where
\begin{equation*}
\gamma ^{s}=\gamma -d\log \left( s\right)
\end{equation*}%
It follows that\ $d\gamma ^{s}=d\gamma $.

One has

\begin{equation}
d\gamma =K\omega _{1}\wedge \omega _{2}.  \label{curvature function}
\end{equation}
This equation is called \emph{the second structure equation of the web, }and%
\emph{\ }the function $K$ is called the \emph{web curvature.}

If we put $d\gamma ^{s}=K^{s}\omega _{1}^{s}\wedge \omega _{2}^{s},$ then $%
K^{s}=s^{-2}K.$ Therefore the curvature function $K$ is a relative invariant
of weight $2.$

In terms of the web function $f,$ one has
\begin{equation}
K=-\frac{1}{f_{x}f_{y}}\left( \log \left( \frac{f_{x}}{f_{y}}\right) \right)
_{xy}  \label{curvature cheres F}
\end{equation}
(cf.\cite{B 55}, \ \S\ 9, or \cite{AS 92}, p. 43).

For the basis vector fields $\partial _{1}$ and $\partial _{2}$, the
structure equations take the form
\begin{equation}
\lbrack \partial _{1},\partial _{2}]=H~(\partial _{2}-\partial _{1}).
\label{vector structure eq}
\end{equation}%
where $\left[ \ ,\right] $ is the commutator of vector fields.

Substituting (\ref{Gamma}) into (\ref{curvature function}), one gets $%
d\gamma =dH\wedge \omega _{1}+\omega _{2}),$and\ from (\ref{differential})
it follows that
\begin{equation}
K=\partial _{1}\left( H\right) -\partial _{2}\left( H\right) .
\label{curvature main}
\end{equation}

\subsection{The Chern Connection}

Let us use the differential $1$-form $\gamma $ to define a connection in the
cotangent bundle $\tau ^{\ast }:T^{\ast }M\rightarrow M$ by the following
covariant differential:
\begin{equation*}
d_{\gamma }:\Lambda ^{1}\left( M\right) \rightarrow \Lambda ^{1}\left(
M\right) \otimes \Lambda ^{1}\left( M\right) ,
\end{equation*}%
where
\begin{eqnarray*}
d_{\gamma }\left( \omega _{1}\right) &\!\!\!\!=&\!\!\!\!-\omega _{1}\otimes
\gamma , \\
d_{\gamma }\left( \omega _{2}\right) &\!\!\!\!=&\!\!\!\!-\omega _{2}\otimes
\gamma ;
\end{eqnarray*}%
and $\otimes $ denotes the tensor product.

In what follows we shall denote by $\Lambda ^{p}\left( M\right) ,$ $p=1,2,$
the modules of smooth differential $p$-forms on $M.$

It is easy to check that the curvature form of the above connection is equal
to $-d\gamma ,$ that is, $d_{\gamma }^{2}:\Lambda ^{1}\left( M\right)
\rightarrow \Lambda ^{1}\left( M\right) \otimes \Lambda ^{2}\left( M\right) $
is the multiplication by $-d\gamma :$%
\begin{equation*}
d_{\gamma }^{2}\left( \omega \right) =-\omega \otimes d\gamma
\end{equation*}%
for any differential form $\omega \in \Lambda ^{1}\left( M\right) .$This
connection is called the \emph{Chern connection} of the web.

It is also easy to check that the Chern connection satisfies the relations
\begin{equation*}
d_{\gamma }\left( \omega _{i}^{s}\right) =-\omega _{i}^{s}\otimes \gamma ^{s}
\end{equation*}%
for $i=1,2,$ and any non-zero smooth function $s.$ The straightforward
computation shows also that $d_{\gamma }$ is a torsion-free connection.

Recall (see, for example, \cite{N 76}, p. 128) that for the covariant
differential $d_{\nabla }:\Lambda ^{1}\left( M\right) \rightarrow \Lambda
^{1}\left( M\right) \otimes \Lambda ^{1}\left( M\right) $ of any
torsion-free connection $\nabla ,$ one has $d_{\nabla }=d_{\gamma }-T,$
where
\begin{equation*}
T:\Lambda ^{1}\left( M\right) \rightarrow S^{2}\left( M\right) \subset
\Lambda ^{1}\left( M\right) \otimes \Lambda ^{1}\left( M\right)
\end{equation*}%
is the \emph{deformation tensor }of the connection, and $S^{2}\left(
M\right) $ is the module of the symmetric $2$-tensors on $M$.

Below we shall use the notation $\nabla _{X}\left( \theta \right) \overset{%
\text{def}}{=}\left( d_{\nabla }\theta \right) \left( X\right) $ for the
covariant derivative of a differential $1$-form $\theta $ along vector field
$X$ with respect to connection $\nabla .$

\begin{proposition}
\label{geodesics proposition}Let $d_{\nabla }:\Lambda ^{1}\left( M\right)
\rightarrow \Lambda ^{1}\left( M\right) \otimes \Lambda ^{1}\left( M\right) $
be the covariant differential of a connection $\nabla $ in the cotangent
bundle of $M.$ Then a foliation $\left\{ \theta =0\right\} $ on $M$ given by
the differential $1$-form $\theta \in \Lambda ^{1}\left( M\right) $ consists
of geodesics of $\nabla $ if and only if
\begin{equation*}
d_{\nabla }\left( \theta \right) =\alpha \otimes \theta +\theta \otimes \beta
\end{equation*}%
for some differential $1$-forms $\alpha ,\beta \in \Omega ^{1}\left(
M\right) .$
\end{proposition}

\begin{proof}
Let $\theta ^{\prime }$ be a differential $1$-form such that $\theta $ and $%
\theta ^{\prime }$ are linearly independent.

Then
\begin{equation*}
d_{\nabla }\left( \theta \right) =\alpha \otimes \theta +\theta \otimes
\beta +h\theta ^{\prime }\otimes \theta ^{\prime }.
\end{equation*}%
Assume that $X$ is a geodesic vector field on $M$ such that $\theta \left(
X\right) =0.$ Then $\nabla _{X}\left( \theta \right) $ must be equal to zero
on $X.$ But
\begin{equation*}
d_{\nabla }\theta \left( X\right) =\beta \left( X\right) \theta +h\theta
^{\prime }\left( X\right) \theta ^{\prime }.
\end{equation*}%
Therefore, $h=0.$
\end{proof}

\begin{corollary}
Foliations $\left\{ \omega _{1}=0\right\} ,\left\{ \omega _{2}=0\right\} ,$
and $\left\{ \omega _{3}=0\right\} $ are geodesic with respect to the Chern
connection.
\end{corollary}

\subsection{Akivis--Goldberg Equations}

The problem of linearization of webs can be reformulated as follows: \emph{%
find a torsion-free flat connection such that the foliations of the web are
geodesic with respect to this connection.}

\begin{proposition}
Let $d_{\nabla }=d_{\gamma }-T:\Lambda ^{1}\left( M\right) \rightarrow
\Lambda ^{1}\left( M\right) \otimes \Lambda ^{1}\left( M\right) $ be the
covariant differential of a torsion-free connection $\nabla $ such that the
foliations $\left\{ \omega _{p}=0\right\} ,\ p=1,2,3,$ are geodesic. Then
\begin{equation}
\begin{array}{ll}
T\left( \omega _{1}\right) =T_{11}^{1}\omega _{1}\otimes \omega
_{1}+T_{12}^{1}\left( \omega _{1}\otimes \omega _{2}+\omega _{2}\otimes
\omega _{1}\right) , &  \\
T\left( \omega _{2}\right) =T_{22}^{2}\omega _{2}\otimes \omega
_{2}+T_{12}^{2}\left( \omega _{1}\otimes \omega _{2}+\omega _{2}\otimes
\omega _{1}\right) \label{defTensor general}, &
\end{array}%
\end{equation}%
where the components of the deformation tensor have the form
\begin{equation}
T_{12}^{2}=\lambda _{1},\ \ T_{12}^{1}=\lambda _{2},\ T_{11}^{1}=2\lambda
_{1}+\mu ,\ T_{22}^{2}=2\lambda _{2}-\mu \   \label{Deformation tensor}
\end{equation}%
for some smooth functions $\lambda _{1},\lambda _{2},$ and $\mu .$
\end{proposition}

\begin{proof}
Due to (\ref{geodesics proposition}) and the requirement that the foliations
$\left\{ \omega _{1}=0\right\} \ $and $\left\{ \omega _{2}=0\right\} \ $are
geodesic, one gets (\ref{defTensor general}). \ The same requirement for the
foliation $\left\{ \omega _{3}=0\right\} $ gives the following relation for
the components of the deformation tensor $T:$
\begin{equation*}
T_{11}^{1}+T_{22}^{2}=2(T_{12}^{1}+T_{12}^{2}),
\end{equation*}
and this implies (\ref{Deformation tensor}).
\end{proof}

\textit{Therefore, in order to linearize the }$3$\textit{-web, one should
find functions }$\lambda _{1},\lambda _{2}\ \ $\textit{and }$\mu $\textit{\
in such a way that the connection corresponding to }$d_{\gamma }-T,$\textit{%
\ where the deformation tensor }$T$\textit{\ has form }(\ref{Deformation
tensor})\textit{, is flat.}

Let us denote by $\nabla _{i}$ the covariant derivatives along $\partial
_{i},i=1,2,$ with respect to the connection $\nabla $ and by
\begin{equation*}
R:\Lambda ^{1}\left( M\right) \rightarrow \Lambda ^{1}\left( M\right)
\end{equation*}%
the curvature tensor .

From the standard formula for the curvature $R\left( X,Y\right) =[\nabla
_{X},\nabla _{Y}]-\nabla _{\lbrack X,Y]}$ (see, for example, \cite{KN 63},
p. 133) and (\ref{vector structure eq}) we find that
\begin{equation*}
R\left( \omega \right) =[\nabla _{1},\nabla _{2}]\left( \omega \right)
+H~(\nabla _{1}-\nabla _{2})\left( \omega \right)
\end{equation*}%
for any $\omega \in \Lambda ^{1}\left( M\right) .$

It follows from the above proposition that for the connection corresponding
to $d_{\gamma }-T$ we get
\begin{eqnarray*}
\nabla _{1}\left( \omega _{1}\right) &\!\!\!\!=&\!\!\!\!-\left( 2\lambda
_{1}+\mu +H\right) ~\omega _{1}-\lambda _{2}~\omega _{2}, \\
\nabla _{1}\left( \omega _{2}\right) &\!\!\!\!=&\!\!\!\!-(\lambda
_{1}+H)~\omega _{2}, \\
\nabla _{2}\left( \omega _{1}\right) &\!\!\!\!=&\!\!\!\!-\left( \lambda
_{2}+H\right) ~\omega _{1}, \\
\nabla _{2}\left( \omega _{2}\right) &\!\!\!\!=&\!\!\!\!-\lambda _{1}~\omega
_{1}-\left( 2\lambda _{2}-\mu +H\right) ~\omega _{2}.
\end{eqnarray*}%
and
\begin{eqnarray*}
R(\omega _{1}) &\!\!\!\!=&\!\!\!\!\left( 2\partial _{2}\left( \lambda
_{1}\right) -\partial _{1}\left( \lambda _{2}\right) +\partial _{2}\left(
\mu \right) -H\left( 2\lambda _{1}-\lambda _{2}+\mu \right) -\lambda
_{1}\lambda _{2}-K\right) ~\omega _{1}+ \\
&&\left( \partial _{2}\left( \lambda _{2}\right) +\lambda _{2}\left(
-H-\lambda _{2}+\mu \right) \right) ~\omega _{2}, \\
R\left( \omega _{2}\right) &\!\!\!\!=&\!\!\!\!(-\partial _{1}\left( \lambda
_{1}\right) +\lambda _{1}\left( H+\lambda _{1}+\mu \right) )~\omega _{1}+ \\
&&\left( \partial _{2}\left( \lambda _{1}\right) -2\partial _{1}\left(
\lambda _{2}\right) +\partial _{1}\left( \mu \right) -H\left( \lambda
_{1}-2\lambda _{2}+\mu \right) +\lambda _{1}\lambda _{2}-K\right) ~\omega
_{2}
\end{eqnarray*}

Therefore, in order to obtain a flat torsion-free connection, components of
the deformation tensor must satisfy the following \emph{Akivis-Goldberg
equations}
\begin{equation}
R\left( \omega _{1}\right) =0,\ R\left( \omega _{2}\right) =0.
\label{Akivis Equation}
\end{equation}
Since $\omega _{1}$\ and $\omega _{2}$\ are linearly independent, equations (%
\ref{Akivis Equation}) imply that
\begin{eqnarray*}
2\partial _{2}\left( \lambda _{1}\right) -\partial _{1}\left( \lambda
_{2}\right) +\partial _{2}\left( \mu \right) -H\left( 2\lambda _{1}-\lambda
_{2}+\mu \right) -\lambda _{1}\lambda _{2}-K &\!\!\!\!=&\!\!\!\!0, \\
\partial _{2}\left( \lambda _{2}\right) +\lambda _{2}\left( -H-\lambda
_{2}+\mu \right) &\!\!\!\!=&\!\!\!\!0, \\
-\partial _{1}\left( \lambda _{1}\right) +\lambda _{1}\left( H+\lambda
_{1}+\mu \right) &\!\!\!\!=&\!\!\!\!0, \\
\partial _{2}\left( \lambda _{1}\right) -2\partial _{1}\left( \lambda
_{2}\right) +\partial _{1}\left( \mu \right) -H\left( \lambda _{1}-2\lambda
_{2}+\mu \right) +\lambda _{1}\lambda _{2}-K &\!\!\!\!=&\!\!\!\!0.
\end{eqnarray*}

Resolving the system with respect to the derivatives of $\lambda _{1}$ and $%
\lambda _{2}$, we obtain the following system of PDEs:
\begin{eqnarray*}
\partial _{1}\left( \lambda _{1}\right) &\!\!\!\!=&\!\!\!\!\lambda
_{1}\left( H+\lambda _{1}+\mu \right) , \\
\partial _{2}\left( \lambda _{1}\right) &\!\!\!\!=&\!\!\!\!\frac{K}{3}%
+H\left( \lambda _{1}+\frac{\mu }{3}\right) +\lambda _{1}\lambda _{2}+\frac{1%
}{3}\partial _{1}\left( \mu \right) -\frac{2}{3}\partial _{2}\left( \mu
\right) , \\
\partial _{1}\left( \lambda _{2}\right) &\!\!\!\!=&\!\!\!\!-\frac{K}{3}%
+H\left( \lambda _{2}-\frac{\mu }{3}\right) +\lambda _{1}\lambda _{2}+\frac{2%
}{3}\partial _{1}\left( \mu \right) -\frac{1}{3}\partial _{2}\left( \mu
\right) , \\
\partial _{2}\left( \lambda _{2}\right) &\!\!\!\!=&\!\!\!\!\lambda
_{2}\left( H+\lambda _{2}-\mu \right) .
\end{eqnarray*}

We shall look at the above system as a system of partial differential
equations with respect to the functions $\lambda _{1}$ and $\lambda _{2}$
provided that $\mu $ is given.

We get the compatibility conditions for this system from structure equations
(\ref{vector structure eq}) for $\lambda _{1}$ and $\lambda _{2}$ presented
in the form
\begin{equation*}
\partial _{1}(\partial _{2}\left( \lambda _{i}\right) )-\partial
_{2}(\partial _{1}\left( \lambda _{i}\right) )+H\left( \partial _{1}\left(
\lambda _{i}\right) -\partial _{2}\left( \lambda _{i}\right) \right) =0,
\end{equation*}%
where $i=1,2.$

After a series of long and straightforward computations, we obtain the
following two compatibility equations:
\begin{equation}
I_{1}\left( \mu \right) =0,\ I_{2}\left( \mu \right) =0,\
\label{Compability Equations}
\end{equation}%
where $I_{1}\left( \mu \right) $ and $I_{2}\left( \mu \right) $\ have the
form
\begin{eqnarray}
&&I_{1}\left( \mu \right) =-\partial _{1}^{2}\left( \mu \right) +2\partial
_{1}\partial _{2}\left( \mu \right) +\left( \mu +H\right) \partial
_{1}\left( \mu \right) -2\left( 2H+\mu \right) \partial _{2}\left( \mu
\right)  \notag \\
&&+H\mu ^{2}+(2H^{2}-\partial _{2}\left( H\right) )\mu -\partial _{1}\left(
K\right) +2HK  \notag
\end{eqnarray}%
and
\begin{eqnarray}
&&I_{2}\left( \mu \right) =-\partial _{2}^{2}\left( \mu \right) +2\partial
_{1}\partial _{2}\left( \mu \right) +2(\mu -H)\partial _{1}\left( \mu
\right) -(H+\mu )\partial _{2}\left( \mu \right) -H\mu ^{2}  \notag \\
&&+\left( 2H^{2}-\partial _{1}\left( H\right) \right) \mu -\partial
_{2}\left( K\right) +2HK.  \notag
\end{eqnarray}

We sum up these results in the following

\begin{theorem}
The Akivis-Goldberg equations as differential equations with respect to the
components $T_{12}^{1}=\lambda _{2}$ and $T_{12}^{2}=\lambda _{1}$ of the
deformation tensor $T$\ are compatible if and only if the component $\mu $
satisfies the following differential equations:
\begin{equation*}
I_{1}\left( \mu \right) =0,\ I_{2}\left( \mu \right) =0.
\end{equation*}%
If the above conditions $(\ref{Compability Equations})$ are valid, then the
system $(\ref{Akivis Equation})$ of \ PDEs is the Frobenius-type system, and
for given values $\lambda _{1}\left( x_{0}\right) $ and $\lambda _{2}\left(
x_{0}\right) $ \ at a point $x_{0}\in M,$ there is (a unique) smooth
solution of the system in some neighborhood of $x_{0}.$
\end{theorem}

It is worthwhile to note the peculiarity of the Akivis-Goldberg system of
differential equations and our presentation of components of the deformation
tensor. This is a non-linear overdetermined system with respect to
components $\lambda _{1},\lambda _{2},\mu $ of the deformation tensor, but
the compatibility conditions in our case depend on $\mu $ only while for
general systems they depend on all components of the deformation tensor.
This gives us a method to find the linearizability conditions in a
constructive way.

\section{Linearizability of 4-Webs}

\subsection{The Basic Invariant of a 4-Web}

A $4$-web $W_{4}$ on $M^{2}$ can be defined by $4$ differential \thinspace $%
1 $-forms $\omega _{1},\omega _{2},\omega _{3},$ and $\omega _{4}$\ such
that any two of them are linearly independent.

We prove the following proposition:

\begin{proposition}
\bigskip The forms $\omega _{1},\omega _{2},\omega _{3},$ and $\omega _{4}$\
can be normalized in such a way that the normalization condition $(\ref%
{normalization equation})$\ holds for the first three of them, and in
addition, the following condition holds for the forms $\omega _{1},\omega
_{2},$ and $\omega _{4:}$\qquad \qquad
\begin{equation}
\omega _{4}+a\omega _{1}+\omega _{2}=0,\   \label{4web normalization}
\end{equation}%
where $a$ is a nonzero function.
\end{proposition}

\bigskip

\begin{proof}
\bigskip In fact, if we take the forms $\omega _{1}$ and $\omega _{2}$ as
cobasis forms of\ $M^{2}$, then the forms $\omega _{3}$ and $\omega _{4}$
are linearly expressed in terms of $\omega _{1}$ and $\omega _{2}:$

\begin{equation*}
\omega _{3}=\alpha \omega _{1}+\beta \omega _{2\,},
\end{equation*}
\begin{equation*}
\omega _{4}=\alpha ^{\prime }\omega _{1}+\beta ^{\prime }\omega _{2\,},
\end{equation*}%
where $\alpha ,\beta ,\alpha ^{^{\prime }},\beta ^{^{\prime }}\neq 0,$%
{\normalsize \ }$\alpha \neq \alpha ^{\prime },\ \alpha \beta ^{^{\prime
}}-\alpha ^{^{\prime }}\beta \neq 0.$

Making the substitution
\begin{equation*}
\omega _{1}\rightarrow -\frac{1}{\alpha }\omega _{1},\ \omega
_{2}\rightarrow \frac{1}{\beta }\omega _{2},\ \omega _{3}\rightarrow -\omega
_{3},\ \omega _{4}\rightarrow -\frac{\beta ^{^{\prime }}}{\beta \mathrm{\ }}%
\omega _{4,}
\end{equation*}%
{\normalsize \ }we get ($\ref{normalization equation}$\ ) and (\ref{4web
normalization}) with $a=\frac{\alpha ^{^{\prime }}\beta }{\beta ^{^{\prime
}}\alpha \mathrm{\ }}$.
\end{proof}

Note that $a\neq 0,1.$ Moreover, the value $a\left( x\right) ,$ $x\in M,$ of
the function $a$ is the cross-ratio of the four tangents to the lines in $%
T_{x}^{\ast }(M^{2})$ \ generated by the covectors $\omega _{1,x},\omega
_{2,x},\omega _{3,x}$, and $\omega _{4,x},$ and therefore is an invariant\
of the $4$-web. The function $a$ is called the \emph{basic invariant }of the
$4$-web (see \cite{G 80} and \cite{G 88}, pp. 302--303).

\subsection{The Expression for $\protect\mu $}

We shall consider a $4$-web $\left\langle \omega _{1},\omega _{2},\omega
_{3},\omega _{4}\right\rangle $ as the $3$-web $\left\langle \omega
_{1},\omega _{2},\omega _{3}\right\rangle $ and an extra foliation given by
form $\omega _{4}$ which satisfies {\normalsize (\ref{4web normalization}).}
Moreover, by the Chern connection, the curvature, etc. that we discussed
above for a $3$-web we shall mean the corresponding constructions for the $3$%
-web $\left\langle \omega _{1},\omega _{2},\omega _{3}\right\rangle .$

\begin{theorem}
Let $\nabla $ be a torsion-free connection in the cotangent bundle $\tau
^{\ast }:T^{\ast }M\rightarrow M$ such that the foliations $\left\{ \omega
_{1}=0\right\} ,\left\{ \omega _{2}=0\right\} ,\left\{ \omega _{3}=0\right\}
,$ and $\left\{ \omega _{4}=0\right\} $ are geodesic for $\nabla .$ Then the
components of the deformation tensor $T$\ have the form $(\ref{Deformation
tensor})$ and
\begin{equation}
\mu =\frac{\partial _{1}\left( a\right) -a\partial _{2}\left( a\right) }{%
a-a^{2}}.  \label{mu-relation}
\end{equation}
\end{theorem}

\begin{proof}
Let $d_{\nabla }=d_{\gamma }-T$ be the covariant differential of the
connection $\nabla $. Then (\ref{4web normalization}) gives
\begin{equation*}
-d_{\nabla }\left( \omega _{4}\right) =\omega _{1}\otimes da-\omega
_{4}\otimes \gamma -aT\left( \omega _{1}\right) -T\left( \omega _{2}\right) .
\end{equation*}%
If $\omega _{4}=0,$ then $\omega _{2}=-a\omega _{1},$ and the right-hand
side takes the form
\begin{equation*}
\left( \partial _{1}\left( a\right) -a\partial _{2}\left( a\right) +\mu
\left( a^{2}-a\right) \right) \omega _{1}\otimes \omega _{1}.
\end{equation*}%
Therefore, this tensor equals zero if and only if equation\ (\ref%
{mu-relation}) holds.
\end{proof}

\bigskip

Formula (16) shows that the quantity $\mu $ occurring in expressions (12) of
the components of the deformation tensor, is expressed in terms of the basic
invariant $a$ and its derivatives. Namely this fact made it possible to
express the linearizability conditions for $4$-webs in terms of $4$th order
jets and solve the linearizability problem for 4-webs without use of
computers.

\section{Differential Invariants of 4-Webs}

For the values of the operators $I_{1}\left( \mu \right) $ and $I_{2}\left(
\mu \right) $ on the function $\mu =(\partial _{1}\left( a\right) -a\partial
_{2}\left( a\right) )/(a-a^{2}),$ we introduce the following operators:
\begin{equation*}
I_{1}^{0}\left( f,a\right) =I_{1}\left( \frac{\partial _{1}\left( a\right)
-a\partial _{2}\left( a\right) }{a-a^{2}}\right)
\end{equation*}%
and%
\begin{equation*}
I_{2}^{0}\left( f,a\right) =I_{2}\left( \frac{\partial _{1}\left( a\right)
-a\partial _{2}\left( a\right) }{a-a^{2}}\right) .
\end{equation*}

These are differential operators of order three in the basic invariant $a$
and of order four in the web function $f.$ If they are equal to zero, then $%
\mu $ satisfies the conditions $I_{1}(\mu )=I_{2}(\mu )=0$, and therefore
the Akivis--Goldberg equations for the $3$-web generated by $\omega
_{1},\omega _{2},$ and $\omega _{3}$ are compatible. They can be solved with
respect to the functions $\lambda _{1}$ and $\lambda _{2},$ and we get
finally the deformation tensor and such a flat connection in which the
leaves of $\omega _{p}=0$ for all $p=1,2,3,4$ are geodesics.

Summarizing we get the following theorem.

\begin{theorem}
The $4$-web $W_{4}$ is linearizable if and only if the conditions $%
I_{1}^{0}\left( f,a\right) =0$ and $I_{2}^{0}\left( f,a\right) =0$ hold.
\end{theorem}

We call the\ quantities $I_{1}^{0}\left( f,a\right) $ and $I_{2}^{0}\left(
f,a\right) $ the\textbf{\ }\emph{basic differential invariants} of the $4$%
-web $W_{4}$.

In order to make the expressions for these invariants more symmetric, we
introduce a second web function for a $4$-web $W_{4}$. Namely, locally one
can find a function $g(x,y)$ such that $\omega _{4}\wedge dg=0,$ or
\begin{equation*}
\omega _{4}=u~dg
\end{equation*}
for some function $u.$ Note that the function $f(x,y)$ defines the $3$%
-subweb of the $4$-web $W_{4}$ formed by the foliations $\left\{ \omega
_{1}=0\right\} ,\left\{ \omega _{2}=0\right\} ,\ $and $\left\{ \omega
_{3}=0\right\} ,$and the function $g(x,y)$ defines the $3$-subweb of the $4$%
-web $W_{4}$ formed by the foliations $\left\{ \omega _{1}=0\right\}
,\left\{ \omega _{2}=0\right\} ,\ $and $\left\{ \omega _{4}=0\right\} .$

It follows from (\ref{4web normalization}) that
\begin{equation*}
ug_{x}=-af_{x},\ ug_{y}=-f_{y}.
\end{equation*}
These two equations imply that
\begin{equation*}
a=\frac{f_{y}g_{x}}{f_{x}g_{y}}
\end{equation*}
and
\begin{equation}
a=\frac{\partial _{1}\left( g\right) }{\partial _{2}\left( g\right) }.
\label{a-invariant}
\end{equation}

Substituting this expression into (\ref{mu-relation}) and the result
obtained into (\ref{Compability Equations}), one gets two differential
invariants $I_{1}\left( f,g\right) $ and $I_{2}\left( f,g\right) $ each of
which is of order three in $f$ and $g.$\bigskip

\subsection{Computation of the Differential Invariants}

\subsubsection{ Calculus of Covariant Derivatives}

Let $d_{\gamma }:\Lambda ^{1}(M)\rightarrow \Lambda ^{1}\left( M\right)
\otimes \Lambda ^{1}\left( M\right) $ be the covariant differential with
respect to the Chern connection.

Denote by $\Theta ^{k}\left( M\right) =\left( \Lambda ^{1}\left( M\right)
\right) ^{\otimes k}$ the module of covariant tensors of order $k.$ Then the
Chern connection induces a covariant differential
\begin{equation*}
d_{\gamma }^{(k)}:\Theta ^{k}\left( M\right) \rightarrow \Theta ^{k+1}\left(
M\right) ,
\end{equation*}%
where
\begin{equation*}
d_{\gamma }^{(k)}:h\theta \longmapsto hd_{\nabla }^{(k)}\left( \theta
\right) +\theta \otimes dh
\end{equation*}%
and $h\in C^{\infty }\left( M\right) \ $and $\theta \in \Theta ^{k}\left(
M\right) .$

If $\theta $ \ has the form $\theta =u\omega _{i_{1}}\otimes \omega
_{i_{2}}\otimes \cdots \otimes \omega _{i_{k}}$ in the basis $\{\omega
_{1},\omega _{2}\},$ where $i_{1},i_{2},...,i_{k}=1,2,$ and $u\in C^{\infty
}\left( M\right) ,$ then
\begin{equation*}
d_{\gamma }^{(k)}\left( \theta \right) =\omega _{i_{1}}\otimes \omega
_{i_{2}}\otimes \cdots \otimes \omega _{i_{k}}\otimes \left( du-ku\gamma
\right) .
\end{equation*}%
We say that $u$ is of weight $k$ and call the form
\begin{equation}
\delta ^{\left( k\right) }\left( u\right) =du-ku\gamma
\label{covariant differential Vadim}
\end{equation}%
the \emph{covariant differential }of $u.$ Decomposing the form $\delta
^{\left( k\right) }\left( u\right) $ in the basis $\{\omega _{1},\omega
_{2}\},$ we obtain
\begin{equation*}
\delta ^{\left( k\right) }\left( u\right) =\delta _{1}^{\left( k\right)
}\left( u\right) ~\omega _{1}+\delta _{2}^{\left( k\right) }\left( u\right)
~\omega _{2},
\end{equation*}%
where
\begin{eqnarray}
\delta _{1}^{\left( k\right) }\left( u\right) &\!\!\!\!=&\!\!\!\!\partial
_{1}\left( u\right) -kHu,  \label{covariant derivatives weight k} \\
\delta _{2}^{\left( k\right) }\left( u\right) &\!\!\!\!=&\!\!\!\!\partial
_{2}\left( u\right) -kHu  \notag
\end{eqnarray}%
are the covariant derivatives of $u$ with respect to the Chern connection.
Note that $\delta _{1}^{\left( k\right) }\left( u\right) $ and $\delta
_{2}^{\left( k\right) }\left( u\right) $ are of weight $k+1.$

\begin{lemma}
For any $s=0,1,...,$ the relation
\begin{equation}
\delta _{2}^{\left( s+1\right) }\circ \delta _{1}^{\left( s\right) }-\delta
_{1}^{\left( s+1\right) }\circ \delta _{2}^{\left( s\right) }=sK
\label{covariant commutator}
\end{equation}%
holds for the commutator.
\end{lemma}

\begin{proof}
We have
\begin{equation*}
\delta _{2}^{\left( s+1\right) }\circ \delta _{1}^{\left( s\right)
}=\partial _{2}\partial _{1}-sH\partial _{2}-\left( s+1\right) H\partial
_{1}+\left( s\left( s+1\right) H^{2}-s\partial _{2}H\right)
\end{equation*}%
and
\begin{equation*}
\delta _{1}^{\left( s+1\right) }\circ \delta _{2}^{\left( s\right)
}=\partial _{1}\partial _{2}-sH\partial _{1}-\left( s+1\right) H\partial
_{2}+\left( s\left( s+1\right) H^{2}-s\partial _{1}H\right) .
\end{equation*}%
The statement follows now from (\ref{curvature main}).
\end{proof}

\subsubsection{Prolongations of the Curvature and the Basic Invariant}

As we have seen, the geometry of a $4$-web is determined by the curvature $K$%
, the basic invariant $a$ and their (covariant) derivatives. In order to
express the invariants $I_{1}$ and $I_{2}$ in terms of $K,a$ and their
covariant derivatives, we need the first covariant derivatives of $K$ and
covariant derivatives of $a$ up to the third order.

We apply (\ref{covariant derivatives weight k}) to $K$ and $a.$

The curvature function $K$ is of weight two. Hence
\begin{eqnarray*}
K_{1} &\!\!\!\!=&\!\!\!\!\delta _{1}^{\left( 2\right) }\left( K\right)
=\partial _{1}\left( K\right) -2HK, \\
K_{2} &\!\!\!\!=&\!\!\!\!\delta _{1}^{\left( 2\right) }\left( K\right)
=\partial _{2}\left( K\right) -2HK.
\end{eqnarray*}

The basic invariant is of weight zero$.$ Hence
\begin{eqnarray*}
a_{1} &\!\!\!\!=&\!\!\!\!\delta _{1}^{\left( 0\right) }\left( a\right)
=\partial _{1}a, \\
a_{2} &\!\!\!\!=&\!\!\!\!\delta _{2}^{\left( 0\right) }\left( a\right)
=\partial _{2}a.
\end{eqnarray*}%
Note that\ (\ref{covariant commutator}) for $s=0$ implies that $\delta
_{2}^{\left( 1\right) }\circ \delta _{1}^{\left( 0\right) }=\delta
_{1}^{\left( 1\right) }\circ \delta _{2}^{\left( 0\right) }.$

Thus, we have the following expressions for the second covariant derivatives
of $a:$%
\begin{eqnarray*}
a_{11} &\!\!\!\!=&\!\!\!\!\delta _{1}^{\left( 1\right) }\circ \delta
_{1}^{\left( 0\right) }(a)=\partial _{1}^{2}a-H\partial _{1}a, \\
a_{12} &\!\!\!\!=&\!\!\!\!a_{21}:=\delta _{2}^{\left( 1\right) }\circ \delta
_{1}^{\left( 0\right) }(a)=\partial _{1}\partial _{2}a-H\partial _{2}a, \\
a_{22} &\!\!\!\!=&\!\!\!\!\delta _{2}^{\left( 1\right) }\circ \delta
_{2}^{\left( 0\right) }(a)=\partial _{2}^{2}a-H\partial _{2}a.
\end{eqnarray*}

Formula (\ref{covariant commutator}) for $s=1$ gives $\delta _{2}^{\left(
2\right) }\circ \delta _{1}^{\left( 1\right) }-\delta _{1}^{\left( 2\right)
}\circ \delta _{2}^{\left( 1\right) }=K.$

Define the third covariant derivatives as follows:
\begin{equation*}
\widetilde{a}_{ijk}=\delta _{k}^{\left( 2\right) }\circ \delta _{j}^{\left(
1\right) }\circ \delta _{i}^{(0)}\left( a\right) .
\end{equation*}%
Note that these expressions are symmetric in $\left( i,j\right) .$ In order
to get symmetry in $\left( i,j,k\right) $ for \emph{all} third covariant
derivatives, we define the \emph{symmetrized third covariant derivatives }$%
a_{ijk}$ as follows:\emph{\ }
\begin{eqnarray*}
a_{111} &\!\!\!\!=&\!\!\!\!\widetilde{a}_{111},a_{222}=\widetilde{a}_{222},
\\
a_{112} &\!\!\!\!=&\!\!\!\!\frac{1}{3}\left( \widetilde{a}_{112}+\widetilde{a%
}_{121}+\widetilde{a}_{211}\right) , \\
a_{122} &\!\!\!\!=&\!\!\!\!\frac{1}{3}\left( \widetilde{a}_{122}+\widetilde{a%
}_{212}+\widetilde{a}_{221}\right) .
\end{eqnarray*}

For them we have the following expressions:
\begin{eqnarray*}
a_{111} &\!\!\!\!=&\!\!\!\!\partial _{1}^{3}a-2H\partial
_{1}^{2}a+(H^{2}-\partial _{1}H)\partial _{1}a, \\
a_{112} &\!\!\!\!=&\!\!\!\!\partial _{1}\partial _{2}\partial
_{1}a-H\partial _{1}^{2}a-2H\partial _{2}\partial _{1}a+\left( 2H^{2}-\frac{%
2\partial _{1}H+\partial _{2}H}{3}\right) \partial _{1}a, \\
a_{122} &\!\!\!\!=&\!\!\!\!\partial _{2}\partial _{1}\partial
_{2}a-H\partial _{2}^{2}a-2H\partial _{1}\partial _{2}a+\left( 2H^{2}-\frac{%
\partial _{1}H+2\partial _{2}H}{3}\right) \partial _{2}a, \\
a_{222} &\!\!\!\!=&\!\!\!\!\partial _{2}^{3}a-2H\partial
_{2}^{2}a+(H^{2}-\partial _{2}H)\partial _{2}a.
\end{eqnarray*}

\subsubsection{ Cartan's Prolongations}

In this section we show the relationship of the above calculus to Cartan's
prolongations of the curvature $K$ and the basic invariant $a$ of a $4$-web $%
W_{4}$.

Since $K$ is a relative invariant of weight two, it satisfies the following
Pfaffian equation:
\begin{equation*}
\delta K=K_{1}\omega _{1}+K_{2}\omega _{2},
\end{equation*}%
where $\delta K=\delta ^{(2)}K=dK-2K\gamma .$

Since $a$ is an absolute invariant, we have
\begin{equation*}
\delta a=a_{1}\omega _{1}+a_{2}\omega _{2},
\end{equation*}%
where $\delta a=\delta ^{(0)}a=da.$

Applying (\ref{covariant differential Vadim}) to $a_{1}\ $and $a_{2},$ we
obtain
\begin{eqnarray*}
\delta a_{1} &\!\!\!\!=&\!\!\!\!a_{11}\omega _{1}+a_{12}\omega _{2}, \\
\delta a_{2} &\!\!\!\!=&\!\!\!\!a_{12}\omega _{1}+a_{22}\omega _{2}
\end{eqnarray*}%
because $a_{12}=a_{21}.$

Here $\delta a_{i}=\delta ^{(1)}a_{i}=da_{i}-a_{i}\gamma ,$ $i=1,2.$

For the covariant differentials of $a_{ij},$ we have
\begin{eqnarray}
\delta a_{11} &\!\!\!\!=&\!\!\!\!\widetilde{a}_{111}\omega _{1}+\widetilde{a}%
_{112}\omega _{2},  \label{2nd covariant differentials} \\
\delta a_{12} &\!\!\!\!=&\!\!\!\!\widetilde{a}_{121}\omega _{1}+\widetilde{a}%
_{122}\omega _{2},  \notag \\
\delta a_{22} &\!\!\!\!=&\!\!\!\!\widetilde{a}_{221}\omega _{1}+\widetilde{a}%
_{222}\omega _{2},  \notag
\end{eqnarray}%
where $\delta a_{ij}=\delta ^{(2)}a_{ij}=da_{ij}-2a_{ij}\gamma .$

Passing to the symmetrized derivatives and using (\ref{covariant commutator}%
) , we find that
\begin{eqnarray*}
\frac{\widetilde{a}_{112}+2\widetilde{a}_{121}}{3} &\!\!\!\!=&\!\!\!%
\!a_{112}, \\
\frac{\widetilde{a}_{112}-\widetilde{a}_{121}}{2} &\!\!\!\!=&\!\!\!\!\frac{K%
}{2}a_{1}.
\end{eqnarray*}%
Therefore,
\begin{equation*}
\widetilde{a}_{112}=a_{112}+\frac{2K}{3}a_{1},
\end{equation*}%
and the first equation in (\ref{2nd covariant differentials}) takes the
following form:
\begin{equation*}
\delta a_{11}=a_{111}\omega _{1}+(a_{112}+\frac{2}{3}a_{1}K)\omega _{2}.
\end{equation*}%
For the second equation of (\ref{2nd covariant differentials}), we have
\begin{equation*}
\widetilde{a}_{121}=a_{112}-\frac{K}{3}a_{1}
\end{equation*}%
and
\begin{equation*}
\delta a_{12}=(a_{112}-\frac{1}{3}a_{1}K)\omega _{1}+\widetilde{a}%
_{122}\omega _{2}.
\end{equation*}%
For the third equation of (\ref{2nd covariant differentials}), we have $%
\widetilde{a}_{122}=\widetilde{a}_{212}$ and
\begin{eqnarray*}
\frac{\widetilde{a}_{221}+2\widetilde{a}_{122}}{3} &\!\!\!\!=&\!\!\!%
\!a_{122}, \\
\frac{\widetilde{a}_{221}-\widetilde{a}_{122}}{2} &\!\!\!\!=&\!\!\!\!-\frac{K%
}{2}a_{2}.
\end{eqnarray*}%
and
\begin{eqnarray*}
\widetilde{a}_{221} &\!\!\!\!=&\!\!\!\!a_{122}-\frac{2}{3}Ka_{2}, \\
\widetilde{a}_{122} &\!\!\!\!=&\!\!\!\!a_{122}+\frac{1}{3}Ka_{2}.
\end{eqnarray*}%
Therefore,
\begin{eqnarray*}
\delta a_{12} &\!\!\!\!=&\!\!\!\!(a_{112}-\frac{1}{3}a_{1}K)\omega
_{1}+(a_{122}+\frac{1}{3}Ka_{2})\omega _{2}, \\
\delta a_{22} &\!\!\!\!=&\!\!\!\!(a_{122}-\frac{2}{3}a_{2}K)\omega
_{1}+a_{222}\omega _{2}.
\end{eqnarray*}

\subsubsection{Differential Invariants in Terms of Covariant Derivatives}

Here we express invariants $I_{1}^{0}\left( f,a\right) $ and $%
I_{2}^{0}\left( f,a\right) $ in terms of the curvature function $K$ , basic
invariant $a$ and their covariant derivatives. To do this, we express the
ordinary derivatives in terms of the covariant derivatines according to the
above formulae. After long computations, we get that the linearizability
conditions $I_{1}^{0}\left( f,a\right) =I_{2}^{0}\left( f,a\right) =0$ are
equivalent to the following two equations:
\begin{eqnarray*}
K_{1} &&\!\!\!\!=\!\!\!\!\frac{1}{a-a^{2}}\Bigl[\frac{1}{3}%
((1-a)a_{1}+aa_{2})K-a_{111}+(2+a)a_{112}-2aa_{122}\Bigr] \\
&&+\frac{1}{(a-a^{2})^{2}}\{[(4-6a)a_{1}+(a^{2}+3a-2)a_{2}]a_{11} \\
&&+[(2a^{2}+7a-6)a_{1}+(2a-3a^{2})a_{2}]a_{12}+[2(a-a^{2})a_{1}-2a^{2}a_{2}]%
\}a_{22} \\
&&+\frac{1}{(a-a^{2})^{3}}[(-6a^{2}+8a-3)(a_{1})^{3}-2a^{3}(a_{2})^{3} \\
&&+(2a^{3}+9a^{2}-15a+6)(a_{1})^{2}a_{2}+(-3a^{3}+6a^{2}-2a)a_{1}(a_{2})^{2}]
\end{eqnarray*}%
and
\begin{eqnarray*}
K_{2} &&\!\!\!\!=\!\!\!\!\frac{1}{a-a^{2}}\Bigl[\frac{1}{3}%
(a_{1}+(a-1)a_{2})K+2a_{112}-(2a+1)a_{122}+aa_{222}\Bigr] \\
&&+\frac{1}{(a-a^{2})^{2}}\{[2a_{1}+(2a-2)a_{2}]a_{11} \\
&&+[(6a-5)a_{1}+(-2a^{2}-3a+2)a_{2}]a_{12}+[(1-a-2a^{2})a_{1}+2a^{2}a_{2}]%
\}a_{22} \\
&&+\frac{1}{(a-a^{2})^{3}}[(4a-2)(a_{1})^{3}+a^{3}(a_{2})^{3} \\
&&+(6a^{2}-12a+5)(a_{1})^{2}a_{2}+(-2a^{3}-3a^{2}+5a-2)a_{1}(a_{2})^{2}].
\end{eqnarray*}

\bigskip

\section{Linearizability of $\boldsymbol{d}$-Webs}

A $d$-web $W_{d}$ on $M^{2}$ is defined by $d$ differential \thinspace $1$%
-forms $\omega _{1},\omega _{2},\omega _{3},...,\omega _{d}$ such that any
two of them are linearly independent. We shall fix the $3$-subweb $%
\left\langle \omega _{1},\omega _{2},\omega _{3}\right\rangle $ and by the
Chern connection, curvature, etc. we shall mean the corresponding
constructions for this $3$-web$.$

For any $4\leq \alpha \leq d,$ we shall consider a $4$-subweb $W_{4}^{\alpha
}$ defined by the forms $\omega _{1},\omega _{2},\omega _{3},\omega _{\alpha
}.$ We denote the basic invariant of this subweb by $a_{\alpha }$ and
continue use the notation $a$ for $a_{4}$. Then
\begin{equation*}
\omega _{\alpha }+a_{\alpha }\,\omega _{1}+\omega _{2}=0.
\end{equation*}

In the same way we used above, we prove the following theorem:

\begin{theorem}
Let $\nabla $ be a torsion-free connection in the cotangent bundle $\tau
^{\ast }:T^{\ast }M\rightarrow M$ such that the foliations $\left\{ \omega
_{1}=0\right\} ,\left\{ \omega _{2}=0\right\} ,\left\{ \omega _{3}=0\right\}
,$ and $\left\{ \omega _{\alpha }=0\right\} $ are $\nabla $-geodesic for all
$\alpha \geq 4.$ Then the components of the deformation tensor $T$\ have
form $(\ref{Deformation tensor})$ and
\begin{equation}
\mu =\frac{\partial _{1}\left( a_{\alpha }\right) -a_{\alpha }\partial
_{2}\left( a_{\alpha }\right) }{a_{\alpha }-a_{\alpha }^{2}}
\label{mu-alpha relation}
\end{equation}
for all $\alpha =4,...,d.$
\end{theorem}

Comparing the expressions for $\mu ,$ we get the following $d-4\;$new
relative invariants of the $d$-web $W_{d}:$
\begin{equation*}
I_{\alpha }=\frac{\partial _{1}\left( a_{\alpha }\right) -a_{\alpha
}\partial _{2}\left( a_{\alpha }\right) }{a_{\alpha }-a_{\alpha }^{2}}-\frac{%
\partial _{1}\left( a\right) -a\partial _{2}\left( a\right) }{a-a^{2}},
\end{equation*}
where $\alpha =5,...,d.$

The web $W_{d}$ can be defined by the functions $f,g_{4}=g,...,g_{d}$ and
\begin{equation*}
a_{\alpha }=\frac{\partial _{1}\left( g_{\alpha }\right) }{\partial
_{2}\left( g_{\alpha }\right) }.
\end{equation*}%
This gives the following expressions for the invariants $I_{\alpha }:$%
\begin{equation*}
I(f,g,g_{\alpha })=I\left( f,g_{\alpha }\right) -I\left( f,g\right) ,
\end{equation*}%
where $\alpha =5,...,d,$ and
\begin{equation*}
I\left( f,p\right) =\frac{(\partial _{1}p)^{2}\,\partial _{2}^{2}p-2\partial
_{1}p\,\partial _{2}p\,\partial _{1}\partial _{2}p+(\partial
_{2}p)^{2}\,\partial _{1}^{2}p}{\partial _{1}p\,\partial _{2}p\,\left(
\partial _{2}p-\partial _{1}p\right) }.
\end{equation*}

Summarizing we get the following theorem:

\begin{theorem}
The $d$-web $W_{d}$ is linearizable if and only if the conditions $%
I_{1}\left( f,g\right) =0$ , $I_{2}\left( f,g\right) =0$ and $I\left(
f,g,g_{5}\right) =0,....,I\left( f,g,g_{d}\right) =0$ hold.
\end{theorem}

\subsection{Method of $\boldsymbol{d}$-Web Linearization}

\subsubsection{4-Webs}

We define a $4$-web $W_{4}$ by two web functions $f$ and $g.$ Then the
procedure for the linearization of such a web can be outlined as follows:

\begin{itemize}
\item[Step 1] Check the linearizability conditions $I_{1}\left( f,g\right)
=0,\;I_{2}\left( f,g\right) =0.$

\item[Step 2] Find the function $\mu $ from (\ref{mu-relation}). Solve the
Akivis-Goldberg equations\ (\ref{Akivis Equation}) with respect to the
functions $\lambda _{1}$ and $\lambda _{2}$. This is the Frobenius-type PDEs
system due to Step 1. Find the components of the deformation tensor $T$ from
(\ref{Deformation tensor}).

\item[Step 3] The connection $\delta _{0}-T$ is flat. Find local coordinates
$x_{1}$ and $x_{2}$ in which the connection coincides with the standard one
on $M^{2}.$ In these coordinates, the leaves of $W_{4}$ are straight lines.
\end{itemize}

\begin{remark}
Step $2$ and Step $3$ can be performed in a constructive way (in
quadratures) if the web under consideration admits a nontrivial symmetry
group. In this case one can find\ the first integrals for the system of
Akivis-Goldberg equations and hence the deformation tensor. If this
deformation tensor also possesses nontrivial symmetries, then the local
coordinates in Step $3$ can be found.
\end{remark}

\subsubsection{$\boldsymbol{d}$-Webs, $\boldsymbol{d>4}$}

We define a $d$-web $W_{d}$ by $d-2$ web functions $f$ and $%
g=g_{4},...,g_{d}.$ Then the procedure for linearization can be outlined as
follows:

\begin{itemize}
\item[Step 1] Check the linearizability conditions $I_{1}\left( f,g\right)
=0,\;I_{2}\left( f,g\right) =0,I(f,g,g_{5})=0,...,I(f,g,g_{d})=0.$

\item[Step 2] Find the function $\mu $ from (\ref{mu-relation}). Solve the
Akivis-Goldberg equations (\ref{Akivis Equation}) with respect to the
functions $\lambda _{1}$ and $\lambda _{2}$. This is the Frobenius-type PDEs
system due to Step 1. Find the components of the deformation tensor $T$ from
(\ref{Deformation tensor}).

\item[Step 3] The connection $\delta _{0}-T$ is flat. Find local coordinates
$x_{1}$ and $x_{2}$ in which the connection coincides with the standard one
on $M^{2}.$ In these coordinates, the leaves of $W_{d}$ are straight lines.
\end{itemize}

\section{ Tests and Examples}

\subsection{Test Notebooks}

Below we give Mathematica codes for testing $4$- and $5$-webs for
linearizability.

The following program computes differential invariants of $d$-webs for $%
d\geq 4$:
\begin{eqnarray*}
&&\mathbf{webInvariants}[fTab\_]:=[\{f,g,X,Y,h,A,I1,I2,J,a,\mu ,d,ans\}, \\
&&f=fTab[[1]];\ d=Length[fTab];\ g[i\_]=fTab[[i]]; \\
&&X[A\_]:=-\frac{D[A,x]}{D[f,x]};\ Y[A\_]:=-\frac{D[A,y]}{D[f,y]};\ h=\frac{%
D[f,x,y]}{D[f,x]\ast D[f,y]}; \\
&&a[i\_]=\frac{D[f,y]\ast D[g[i],x]}{D[f,x]\ast D[g[i],y]};\ \ \nu \lbrack
i\_]:=\frac{X[a[i]]-a[i]\ast Y[a[i]]}{a[i]^{2}-a[i]};\ \mu =\nu \lbrack 2];
\\
&&I1=X[X[\mu ]]-2\ast X[Y[\mu ]]+(\mu -h)\ast X[h]+(4\ast h-2\ast \mu )\ast
Y[\mu ]+ \\
&&h\ast \mu ^{2}-(2\ast h^{2}-Y[h])\ast \mu -X[X[h]]+X[Y[h]]+2\ast h\ast X[h]
\\
&&-2\ast h\ast Y[h]//\mathbf{Simplify;} \\
&&I2=X[Y[\mu ]]-2\ast X[Y[\mu ]]+(2\ast \mu +2\ast h)\ast X[h]+(h-\mu )\ast
Y[\mu ]- \\
&&h\ast \mu ^{2}-(2\ast h^{2}-X[h])\ast \mu +Y[Y[h]]-Y[X[h]]+2\ast h\ast X[h]
\\
&&-2\ast h\ast Y[h]//\mathbf{Simplify;} \\
&&J[i\_]:=(\mu -\nu \lbrack i])//\mathbf{Simplify;} \\
&&ans=\{I1,I2,Table[J[i],\{i,3,d\}]\}\ \ ]
\end{eqnarray*}

The following\ program tests $4$-webs for the linearizability:%
\begin{eqnarray*}
&&\mathbf{LinTest4Web}[f\_,g\_]:=\mathbf{Module}[ \\
&&\{X,Y,h,A,I1,I2,a,\mu ,Z,ans\}, \\
&&X[A\_]:=-\frac{D[A,x]}{D[f,x]};\ Y[A\_]:=-\frac{D[A,y]}{D[f,y]};\ h= \\
&&\frac{D[f,x,y]}{D[f,x]\ast D[f,y]}; \\
&&a=\frac{D[f,y]\ast D[g,x]}{D[f,x]\ast D[g,y]};\ \mu =\frac{X[a]-a\ast Y[a]%
}{a^{2}-a}; \\
&&I1=X[X[\mu ]]-2\ast X[Y[\mu ]]+(\mu -h)\ast X[h]+(4\ast h-2\ast \mu )\ast
Y[\mu ]+ \\
&&h\ast \mu ^{2}-(2\ast h^{2}-Y[h])\ast \mu -X[X[h]]+X[Y[h]]+2\ast h\ast X[h]
\\
&&-2\ast h\ast Y[h]//\mathbf{Simplify;} \\
&&I2=X[Y[\mu ]]-2\ast X[Y[\mu ]]+(2\ast \mu +2\ast h)\ast X[h]+(h-\mu )\ast
Y[\mu ]- \\
&&h\ast \mu ^{2}-(2\ast h^{2}-X[h])\ast \mu +Y[Y[h]]-Y[X[h]]+2\ast h\ast X[h]
\\
&&-2\ast h\ast Y[h]//\mathbf{Simplify;} \\
&&Z=If[I1===0\&\&I2===0,"YES","NO"]; \\
&&ans=Z\ \ ]
\end{eqnarray*}

Finally we give the code which tests $d$-webs.

\begin{eqnarray*}
&&\mathbf{LindTestdWeb}[fun\_]:=\mathbf{Module}[\{f,g,X,Y,h,d,I1,I2,J,a,\nu
,\mu ,Z,ans\}, \\
&&\ f=fun[[1]];d=Length[fun];g[i\_]:=fun[[i]]; \\
&&X[A\_]:=-\frac{D[A,x]}{D[f,x]};\ \ Y[A\_]:=-\frac{D[A,y]}{D[f,y]};\ \ h=%
\frac{D[f,x,y]}{D[f,x]\ast D[f,y]}; \\
&&\ a[i\_]:=\frac{D[f,y]\ast D[g[i],x]}{D[f,x]\ast D[g[i],y]};\ \ \ \nu
\lbrack i\_]:=\frac{X[a[i]]-a[i]\ast Y[a[i]]}{a[i]\symbol{94}2-a[i]};\ \mu
=\nu \lbrack 2]; \\
&&I1=X[X[\mu ]]-2\ast X[Y[\mu ]]+(\mu -h)\ast X[h]+(4\ast h-2\ast \mu )\ast
Y[\mu ]+ \\
&&h\ast \mu ^{2}-(2\ast h^{2}-Y[h])\ast \mu -X[X[h]]+X[Y[h]]+2\ast h\ast
X[h]- \\
&&2\ast h\ast Y[h]//\mathbf{Simplify;} \\
&&I2=X[Y[\mu ]]-2\ast X[Y[\mu ]]+(2\ast \mu +2\ast h)\ast X[h]+(h-\mu )\ast
Y[\mu ]- \\
&&h\ast \mu ^{2}-(2\ast h^{2}-X[h])\ast \mu +Y[Y[h]]-Y[X[h]]+2\ast h\ast
X[h]- \\
&&2\ast h\ast Y[h]//\mathbf{Simplify;} \\
&&J[i\_]:=(\mu -\nu \lbrack i])//\mathbf{Simplify;} \\
&&Z=\mathbf{If}[I1===0\&\&I2===0\&\&\  \\
&&\mathbf{Table[}J[i],\{i,3,d\}]===\mathbf{Table[}0,\{i,3,d\}],"YES","NO"];
\\
&&ans=Z\ \ ]
\end{eqnarray*}

In the last test $fun$ \ is a collection $\{f_{1},...,f_{d-2}\}$ of
functions determining the $d$-web.

Results of the tests are ''YES''\ or ''NO''\ depending on the
linearizability of the web. Note that the computer testing gives the same
results if in each example we replace the functions $f(x,y)$ and $g(x,y)$ by
the functions $f(p(x),q(y))$ and $g(p(x),q(y))$, where $p(x)$ \ and $q(y)$\
are arbitrary smooth functions of $x$ and $y,$ respectively (i.e., if we
consider equivalent webs).

\subsection{Examples}

\begin{enumerate}
\item $\mathbf{LinTest4Web}[x/y,x+y]="YES"$

This is the $\mathbf{4}$-web whose $\mathbf{3}$rd foliation consists of
straight lines of the pencil with center at the origin, and the $\mathbf{4}$%
th foliation consists of parallel straight lines forming the angle $\mathbf{%
135}$ degrees with positive direction of the axis $Ox$, i.e., this $\mathbf{4%
}$-web is linear, and the test is just for demonstration that it is working.

\item $\mathbf{LinTest4Web}[x/y,(1-y)/(1-x)]="YES"$

In this case the $\mathbf{3}$rd and $\mathbf{4}$th foliations are straight
lines of two pencils with their vertices at $(0,0)$ and $(1,1)$. This $%
\mathbf{4}$-web is also linear, and the test is just for demonstration that
it is working.

\item $\mathbf{LinTest4Web}[x+\sqrt{x^{2}-y},x+y]="YES"$

In this case the curves of the $\mathbf{3}$rd foliation are tangent to the
parabola $y=x^{2}$, and the $\mathbf{4}$th foliation consists of parallel
straight lines forming the angle $\mathbf{135}$ degrees with positive
direction of the axis $Ox$, i.e., this $\mathbf{4}$-web is linear. But here
it is not obvious, that the $\mathbf{3}$rd foliation consists of straight
lines.

\item $\mathbf{LinTest4Web}[x+\sqrt{x^{2}-y},y+\sqrt{y^{2}-x}]="YES"$

Here the curves of the $\mathbf{3}$rd foliation are tangent to the parabola $%
y=x^{2}$, and the curves of the $\mathbf{4}$th foliation are tangent to the
parabola $x=y^{2}$, \ i.e., this $\mathbf{4}$-web is linear.

\item $\mathbf{LinTest4Web}[x/y,(x+y)\ast Exp[-x]]="NO"$

This is the $\mathbf{4}$-web whose $\mathbf{3}$rd foliation consists of
straight lines of the pencil with center at the origin, and the $\mathbf{4}$%
-subweb defined by the $\mathbf{4}$th foliation and the coordinate lines is
parallelizable. The $\mathbf{4}$-web in this example is not linearizable,
although two of its $\mathbf{3}$-subwebs are linearizable.

\item $\mathbf{LinTest4Web}[x/y,x^{n}+y^{n}]="YES"$

This web is equivalent to the $\mathbf{4}$-web of the $\mathbf{1}$st
example. This web is not linear but it is linearizable.

\item $\mathbf{LinTestdWeb}[\{y/x,(1-y)/(1-x),(x-xy)/(y-xy)\}]="NO"$

This is the famous $\mathbf{5}$-web constructed by Bol (see \cite{B 55}, \S\
46 and \cite{BB 38}, \S 12 and \S 31). The web consists of $\mathbf{4}$
pencils of straight lines (the first two are the pencils of parallel
coordinate lines, and the $\mathbf{3}$rd and the $\mathbf{4}$th are the
pencils with centers at $(0,0)$ and $(1,1)$), and a foliation of conics
passing through $\mathbf{4}$ centers of the $\mathbf{4}$ pencils. Bol
constructed this example to show that there exists a $\mathbf{5}$-web of
maximum rank $\mathbf{6}$ which is not linearizable. Bol gave an indirect
proof that this $\mathbf{5}$-web is not linearizable. Our test gives the
direct proof of this fact.

\item $\mathbf{LinTest4Web}[y/x,(x-xy)/(y-xy)]="YES"$

This is a $\mathbf{4}$-subweb of the Bol $\mathbf{5}$-web considered in the
previous example. It is formed by $\mathbf{3}$ pencils of straight lines and
the same foliation of conics. It appeared that this $\mathbf{4}$-web is
linearizable while the Bol $\mathbf{5}$-web is not linearizable. Note that
we can prove the linearizability of this $\mathbf{4}$-web using the
quadratic transformation $x=1/x^{\ast },y=1/y^{\ast }$ suggested by Blaschke
in \cite{B 55}, \S 46.

\item $\mathbf{LinTestdWeb}%
[\{x/y,(1-y)/(1-x),(x-xy)/(y-xy),xy,(x-xy)/(x-1),(1-y)/(xy-y),x(1-y)^{2}/y(1-x)^{2}\}]="NO"
$

This is the Spence--Kummer $\mathbf{9}$-web constructed by Pirio and Robert
(see \cite{P 02a}, \cite{P 02b} and \cite{R 02}). This web consists of $%
\mathbf{4}$ pencils of straight lines described in Example 7, \ $\mathbf{4}$
foliations of conics and a foliation of cubics passing through $\mathbf{4}$
centers of the $\mathbf{4}$ pencils. Pirio and Robert constructed this
example and other examples of $d$-webs, $d=6,7,8,$ to show that there exist
nonlinearizable webs of maximum rank different from the Bol $\mathbf{5}$-web
considered in Example $7$. They proved that all their $d$-webs are not
linearizable. Our test gives the direct proof of this fact for the
Spence--Kummer $\mathbf{9}$-web (and all other $d$-webs constructed in \cite%
{P 02a}, \cite{P 02b} and \cite{R 02}).
\end{enumerate}

\bigskip

\bigskip

\emph{Authors' addresses:}

Department of Mathematics, Jerusalem College of Technology-Machon Lev,
Havaad Haleumi St., POB 16031, Jerusalem 91160, Israel; akivis@mail.jct.ac.il

Deparment of Mathematical Sciences, New Jersey Institute of Technology,
University Heights, Newark, NJ 07102, USA; vlgold@oak.njit.edu

Department of Mathematics, The University of Tromso, N9037, Tromso, Norway;
lychagin@math.uit.no

\end{document}